\documentclass[12pt]{amsart}
\usepackage{hyperref}
\usepackage{enumerate}

\begin{document}
\numberwithin{equation}{section}

\def\1#1{\overline{#1}}
\def\2#1{\widetilde{#1}}
\def\3#1{\widehat{#1}}
\def\4#1{\mathbb{#1}}
\def\5#1{\frak{#1}}
\def\6#1{{\mathcal{#1}}}

\def\C{{\4C}}
\def\R{{\4R}}
\def\N{{\4N}}
\def\Z{{\4Z}}

\title[q-effectiveness]
{q-effectiveness for holomorphic subelliptic multipliers
	\smallskip
	\smallskip
	{
		\sl
		D\MakeLowercase
		{edicated to} 
		P\MakeLowercase{rofessor} J.J.~K\MakeLowercase{ohn
		on the occasion of his birthday
		}
	}
}
\author[S.-Y. Kim \& D. Zaitsev]{Sung-Yeon Kim and Dmitri Zaitsev}
\address{S.-Y. Kim: Center for Complex Geometry, Institute for Basic Science,
55 Expo-ro Yuseong-gu Daejeon 34126, South Korea}
\email{sykim8787@ibs.re.kr}
\address{D. Zaitsev: School of Mathematics, Trinity College Dublin, Dublin 2, Ireland}
\email{zaitsev@maths.tcd.ie}
\maketitle

\def\Label#1{\label{#1}}


\def\cn{{\C^n}}
\def\cnn{{\C^{n'}}}
\def\ocn{\2{\C^n}}
\def\ocnn{\2{\C^{n'}}}


\def\dist{{\rm dist}}
\def\const{{\rm const}}
\def\rk{{\rm rank\,}}
\def\id{{\sf id}}
\def\aut{{\sf aut}}
\def\Aut{{\sf Aut}}
\def\CR{{\rm CR}}
\def\GL{{\sf GL}}
\def\Re{{\sf Re}\,}
\def\Im{{\sf Im}\,}
\def\span{\text{\rm span}}
\def\mult{\text{\rm mult\,}}
\def\reg{\text{\rm reg\,}}
\def\ord{\text{\rm ord\,}}
\def\hot{\text{\rm HOT\,}}

\def\codim{{\rm codim}}
\def\crd{\dim_{{\rm CR}}}
\def\crc{{\rm codim_{CR}}}

\def\eps{\varepsilon}
\def\d{\partial}
\def\a{\alpha}
\def\b{\beta}
\def\g{\gamma}
\def\G{\Gamma}
\def\D{\Delta}
\def\Om{\Omega}
\def\k{\kappa}
\def\l{\lambda}
\def\L{\Lambda}
\def\z{{\bar z}}
\def\w{{\bar w}}
\def\Z{{\1Z}}
\def\t{\tau}
\def\th{\theta}

\emergencystretch15pt
\frenchspacing

\newtheorem{Thm}{Theorem}[section]
\newtheorem{Cor}[Thm]{Corollary}
\newtheorem{Pro}[Thm]{Proposition}
\newtheorem{Lem}[Thm]{Lemma}
\newtheorem{Prob}[Thm]{Problem}
\newtheorem{Conj}[Thm]{Conjecture}

\theoremstyle{definition}\newtheorem{Def}[Thm]{Definition}

\theoremstyle{remark}
\newtheorem{Rem}[Thm]{Remark}
\newtheorem{Exa}[Thm]{Example}
\newtheorem{Exs}[Thm]{Examples}

\def\bl{\begin{Lem}}
\def\el{\end{Lem}}
\def\bp{\begin{Pro}}
\def\ep{\end{Pro}}
\def\bt{\begin{Thm}}
\def\et{\end{Thm}}
\def\bc{\begin{Cor}}
\def\ec{\end{Cor}}
\def\bd{\begin{Def}}
\def\ed{\end{Def}}
\def\br{\begin{Rem}}
\def\er{\end{Rem}}
\def\be{\begin{Exa}}
\def\ee{\end{Exa}}
\def\bpf{\begin{proof}}
\def\epf{\end{proof}}
\def\ben{\begin{enumerate}}
\def\een{\end{enumerate}}
\def\beq{\begin{equation}}
\def\eeq{\end{equation}}

\begin{abstract}
We provide a solution to the effectiveness problem
in Kohn's algorithm for generating holomorphic subelliptic multipliers
for $(0,q)$ forms for arbitrary $q$.
%
As application, we obtain subelliptic estimates
for $(0,q)$ forms with effectively controlled
order $\eps>0$ (the Sobolev exponent)
for domains given by sums of squares of holomorphic functions
(J.J. Kohn called them ``special domains'' in \cite{K79}).
These domains 
are of particular interest
due to their relation with 
complex and algebraic geometry.
Our methods include
{\em triangular resolutions} introduced
by the authors in \cite{KZ20}.
%


\end{abstract}

\section{Introduction}

In his celebrated paper \cite{K79},
J.J. Kohn
invented 
a purely algebraic approach
to subelliptic estimates
to the $\bar\d$ problem,
based on {\em generating multiplier ideals}
that, quoting Y.-T.~Siu~\cite{S17},
``measure location and extent of failure of subelliptic estimates'':

\bd
Let $\Om\subset\cn$ be a domain and
$p\in\d\Om$ a boundary point.
\ben
\item
\cite[Definition 1.11]{K79}
A {\em subelliptic estimate} of order $\eps>0$
for $(0,q)$ forms
 is said to hold at $p$
if there exist an open neighborhood
$U$ of $p$  and $C>0$ satisfying
\[
	\|u\|_\eps^2
	\le
	C(
		\|\bar\d u\|^2
		+ \|\bar\d^* u\|^2
		+\|u\|^2
	)
\]
for all $(0,q)$ forms $u$ 
with compact support in $U\cap\1\Om$ 
which belong to the domain of the adjoint $\bar\d^*$.
Here
$\|\cdot\|_\eps$ 
and 
$\|\cdot\|$
are
respectively the tangential Sobolev norm
of the (fractional) order $\eps$
and
the standard $L^2$ norm on $\Omega$.
\item
\cite[Definition 4.2]{K79}
A {\em subelliptic multiplier} of order $\eps>0$ at $p$ for $(0,q)$ forms,
called here briefly a ``{\em $q$-multiplier}'',
is a germ $f$ of a smooth function at $p$,
for which there is a representative in a neighborhood $U$ of $p$,
also denoted by $f$ and $C>0$ satisfying
\[
\Label{subel}
	\|fu\|_\eps^2
	\le
	C(
		\|\bar\d u\|^2
		+ \|\bar\d^* u\|^2
		+\|u\|^2
	)
\]
for all $(0,q)$ forms $u$ as in the subelliptic estimate.
\een
\ed

In particular, a subelliptic estimate of order $\eps>0$ holds at $p$ if and only if
$f=1$ is a {\em $q$-multiplier} of order $\eps>0$ at $p$.
Although multipliers are {\em defined} in terms of the above 
{\em a priori estimate},
J.J.~Kohn discovered in \cite{K79}
purely algebraic procedures 
generating multipliers starting from the defining equation of $\d\Om$.
Based on these procedures, J.J.~Kohn 
proved for bounded domains with
{\em real-analytic} boundary of finite D'Angelo type,
that the trivial multiplier $f=1$ can be generated
by a finite sequence of these procedures.
\cite[Theorem~1.19]{K79}


On the other hand, 
the more general case of
{\em smooth boundary}
remains open (also formulated by Y.\nobreakdash-T.~Siu~\cite[\S2]{S17}):

\begin{Conj}[Kohn's conjecture]
For a bounded pseudoconvex domain 
with smooth boundary
of finite type in $\C^n$,
the trivial multiplier $f=1$
can be generated
by a finite sequence of Kohn's procedures.
\end{Conj}

The stronger {\em Effective Kohn's conjecture}
with the additional control of the order $\eps$
remains open even for real-analytic boundaries:


\begin{Conj}[Effective Kohn's conjecture]
Kohn's conjecture holds under the same assumptions
with an additional effective estimate of the 
{\em order of subellipticity} of the multiplier $f=1$
as function of the finite type and the dimension $n$.
\end{Conj}

The effective conjecture is known for $n=2$, see \cite[\S8]{K79},
where it is based on fundamental results
by H\"ormander \cite{H65} and Rothschild-Stein \cite{RS76}.
In higher dimension, the situation is much less understood,
in fact, examples of \cite{Heier} (\S1.1 in the preprint version)
and \cite[Proposition 4.4]{CD10} in dimension $3$
illustrate a lack of such control,
see also \cite[\S4.1]{S17} for
a detailed explanation of this important phenomenon.

When $n>2$, only the case $q=1$ has been previously considered.
To tackle the effectiveness,
Siu \cite{S10, S17} 
 introduced algebraic geometric techniques
to obtain the effectiveness 
 in the important case of 
{\em special domains}
(see Definition~\ref{spec} below)
 of finite type
in dimension $3$,
with further indications of how to proceed in
the more general cases of
special domains in higher dimension,
and outlining a program to treat 
the more general real-analytic and smooth cases.
A different effective procedure
in Kohn's algorithm 
was given by D'Angelo~\cite{D95}
and Catlin-D'Angelo \cite[Section~5]{CD10}
for special domains given by 
so-called {\em triangular systems}
of holomorphic functions.
In \cite{N14}  A.C.~Nicoara
proposed a construction for the termination of the Kohn algorithm in the real-analytic case with an indication of the ingredients needed for the effectivity.
More recently, the authors of this article 
established another effective procedure
for special domains in dimension $n=3$ \cite{KZ18}
and arbitrary $n$ in \cite{KZ20}
(Y.-T. Siu also told us about his unpublished proof in this case).

The reader is referred to
\cite{D95, DK99, S01, S02, K04, S05, S07, Ch06, S09, CD10, S10, S17,Fa20}
for more extensive details on subelliptic multipliers
and Siu's accounts \cite{S07, S09, S17}
on their broad role and relation with 
other multipliers arising in complex and algebraic geometry.
See also 
\cite{K00, K02, K04, K05, Ce08, St08, CS09, Ba15, BPZ15, CZ17, S17}
for multipliers in more general settings.

\subsection{Main results}

The goal of this note is to provide a solution to the effectiveness problem
in Kohn's algorithm for holomorphic subelliptic multipliers
for 
{\em $(0,q)$ forms for arbitrary $1\le q\le n$}.
We first 
recall 
{\em Kohn's
multiplier generating procedures} 
for holomorphic multipliers
\cite[\S7]{K79}
that
can be described algebraically 
starting with an abstract 
{\em initial set} of germs:

%



\bd[Holomorphic Kohn's procedures]\Label{k-hol}
For an arbitrary initial subset $S$
in the set $\6O_{n,p}$
of holomorphic function germs in $\C^n$ at a point $p\in \C^n$
and integer $1\le q\le n$,
the {\em holomorphic Kohn's ($q$-)procedures}
consist of:
\begin{enumerate}
\item[(P1)]
	for $0<\eps\le 1/2$ and $f_1, \ldots, f_{n-q+1}$ either in $S$
	or multipliers of order $\ge\eps$,
	it follows that
	the partial Jacobian $(n-q+1)\times (n-q+1)$ minors
	$$
		\frac
		{\partial( f_1,\ldots, f_{n-q+1})}
		{\partial (z_{j_1},\ldots, z_{j_{n-q+1}})},
		\quad
		1\le j_1<\ldots< j_{n-q+1}\le n,
	$$
	are multipliers of order $\ge \eps/2$;
	
\item[(P2)] 
	for $0<\eps<1$, $k,r\in\N_{\ge 1}$,  
	$f_1, \ldots, f_k$ multipliers or order $\ge\eps$,
	and
	$g$ a holomorphic function (germ) with $g^r\in (f_1, \ldots, f_k)$,
	it follows that $g$ is a multiplier of order $\ge \eps/r$.
\end{enumerate}
\ed


Rather than directly using the finite type,
we control
the order of subellpticity of $q$-multipliers
 in terms of the {\em $q$-multiplicity}
 defined as follows:

\bd
The {\em $q$-multiplicity} $\mult_q (I)$ of 
an ideal $I\subset\6O_{n,p}$ 
in the ring $\6O_{n,p}$ of germs at $p$ 
of holomorphic functions in $\cn$
is the minimum of the dimension of the quotient space
$$
	\mult_q(I):=\min \dim
	\6 O_{n,p}/( I+(L_1,\ldots, L_{q-1})),
$$
 where the minimum is taken over all choice of $(q-1)$
 affine linear functions $L_1,\ldots, L_{q-1}$ vanishing at $p$.
By the {\em $q$-multiplicity} of
a subset $S\subset \6O_{n,p}$
we mean the $q$-multiplicity of the ideal
generated by the set
$\{ f-f(p) : f\in S\}$.
\ed

Note that $q$-multiplicity is in fact a biholomorphic invariant
(\S\ref{defs}).
We formulate our first
 result 
purely in terms of
Kohn's procedures (P1) and (P2):

\bt\Label{p12}
For every number $\nu>1$ and initial subset $S\subset \6O_{n,p}$ of finite $q$-multiplicity $\le \nu$,
there exists an effectively computable sequence 
$f_1,\ldots, f_m\in \6O_{n,p}$, where $f_m=1$ and each 
$f_j$ is either in $S$ or is obtained by applying 
to $S$ or multipliers from $\{f_1,\ldots,f_{j-1}\}$
one of the Kohn's procedures (P1) or (P2).
Furthermore, the number of steps and the root orders in (P2)
are effectively bounded by functions depending only on $(n,q,\nu)$.
\et

As the first application, 
we obtain
the effectiveness 
for the so-called special domains
 \cite[\S7]{K79}, \cite[\S2.8]{S17}:

\bd\Label{spec}
A {\em special domain} in $\C^{n+1}$
is one defined locally near each boundary point $p$
by
\beq\Label{om}
	\Re(z_{n+1})+\sum_{j=1}^N |F_j(z_1,\ldots, z_n)|^2<0
\eeq
where $F_1,\ldots,F_N$ are holomorphic functions
in a neighborhood of $p$.
By the {\em $q$-multiplicity} of domain \eqref{om} at $p$ 
we mean the $q$-multiplicity of the set
$S=\{F_1,\ldots,F_N\}$.
\ed
 


As an immediate consequence of Kohn's theory and
Theorem~\ref{p12}
applied to $S$ as in Definition~\ref{spec}, we obtain:

\bc\label{e-est-1}
There exists a positive function 
$\eps\colon \N_{>0} \times \N_{>0}\times \N_{>0} \to \R_{>0}$
such that for integers $\nu,n,q\in \N_{>0}$ and any domain 
\eqref{om}
of finite $q$-multiplicity 
$\le \nu$ at a boundary point $p$,
a subelliptic estimate
for $(0,q)$ forms holds at $p$
with effectively bounded
order of subellipticity $\ge \eps(n,q,\nu)$.
\ec

\br
Since the $q$-multiplicity of \eqref{om} is
$ \le (T/2)^{n-q+1}$
where $T$ is the D'Angelo 
{\em $q$-finite type} of \eqref{om} at $p$
by a result of D'Angelo 
\cite[Theorem~2.7]{D82},
an effective bound in terms of the type
can be obtained by
substituting $(T/2)^{n-q+1}$ for $\nu$
in Theorem \ref{e-est-1}.
See also \cite{BS92,BHR96,FIK96,FLZ14,MM17,D17,Fa19, Fa20,HY19,Z19}
for relations of the finite type with other invariants.
\er

\subsection{Triangular resolutions and effective meta-procedures}
In this section we introduce our main tools.
Recall that the crucial lack of effectiveness in  (P2) (see Definition~\ref{k-hol}) is due to the
fact that the order of
subellipticity of
 the generated multiplier depends
  on the root order that can be arbitrarily large in general.

To quantify this phenomenon, we call a procedure
{\em effective} if the order of the new multiplier
can be effectively estimated in terms of 
a quantity associated to the data that we call {\em a complexity}.
We don't seek complexities of individual multipliers 
but rather of their {\em finite tuples and tuples of their ideals,
or more precisely, their filtrations}.
More specifically,
we shall
use
the notion of triangular resolution 
that we defined in the earlier paper \cite{KZ20}:


\bd[\cite{KZ20}]
A {\em triangular resolution of length $k\ge 1$}
and multi-order $(\mu_1,\ldots,\mu_k)\in \N^k$
of a pair $(\G, \6I)$,
where 
$\G\colon (\cn,0)\to (\C^{n},0)$
is a holomorphic map germ 
and
$I_1\subset \ldots \subset I_k \subset \6O_{n,0}$
 a filtration $\6I$ of ideals,
is 
a system of holomorphic function germs
$(h_1,\ldots,h_k)$
satisfying
$$
	h_j = h_j(w_j, \ldots, w_{n}),
	\quad
	h_j\circ \G \in I_j,
	\quad
	\ord_{w_j} h_j = \mu_j,
	\quad
	1\le j\le k.
$$
\ed


%

Our proof of the results from previous section
is based on the following effective meta-procedures
involving triangular resolutions:

\bt\Label{meta}
For integers $n\ge1$, $1\le q \le n$, $0\le k\le n-q$ and $\mu\ge1$, the following hold:
\begin{enumerate}
\item[(MP1)] (Selection of a partial Jacobian).
For any collections of germs
$$
	f=(f_1,\ldots,f_k)\in (\6O_{n,0})^k,
		\quad 
	\psi=(\psi_{k+1}, \ldots, \psi_{n-q+1})\in (\6O_{n,0})^{n-k-q+1},
		\quad
	\mult_q(f,\psi) \le \mu,
$$
there exist linear changes 
of the coordinates $z\in \C^n$
and 
of the components of $\psi$ in $\C^{n-k-q+1}$ 
such that for
the partial Jacobian determinant
$$
	J:= \frac{\d(\psi_{k+1},\ldots,\psi_{n-q+1})}{\d(z_{k+1},\ldots,z_{n-q+1})},
$$
the $q$-multiplicity
$$
	\mult_q
(f,J,\psi_{k+2}, \ldots, \psi_{n-q+1})
$$
is effectively bounded by a function depending only on 
$(n,q,\mu)$.

\item[(MP2)] (Selection of a triangular resolution).
For any collections of germs
$$
	f=(f_1,\ldots,f_k)\in (\6O_{n,0})^{k},
		\quad 
	\psi=(\psi_{1}, \ldots, \psi_{n-q+1})\in (\6O_{n,0})^{n-q+1},
$$
with
$$
	\mult_q(f_1,\ldots,f_j, \psi_{j+1}, \ldots, \psi_{n-q+1})\le \mu
	\quad
	\text{ for all }
	\;
	0\le j\le k,
$$
there exist a germ of a holomorphic map
$$
	\Gamma_\psi:=(\psi, L_{n-q+2},\ldots,L_n)\colon (\mathbb C^n,0)\to (\mathbb C^n,0),
	\quad
	L_j \text{ are linear functions},
$$
such that 
$$
	\mult_q(\psi)=\mult(\Gamma_\psi)
$$
and a triangular resolution $h=(h_1,\ldots,h_k)$ of $(\Gamma_\psi,\6I)$,
where $\6I$ is the filtration 
$$
	(f_1)\subset (f_1, f_2) \ldots \subset (f_1,\ldots, f_k),
$$
such that orders $\ord_{w_j}h_j$ are 
effectively bounded by functions depending only on 
$(n,q,\mu)$.

\item[(MP3)] (Jacobian extension in a triangular resolution).
For any 
$$
\G=(\phi,\psi, z_{n-q+2},\ldots, z_n),
\quad (\phi,\psi)\in (\6O_{n,0})^{k}\times (\6O_{n,0})^{n-k-q+1}
$$
and filtration $\6I$ of ideals 
$I_1\subset \ldots \subset I_{k+1}\subset \6O_{n,0}$
satisfying 
$$
	I_{k+1}\subset I_k + (J),
$$ 
where $J$ is the Jacobian determinant of $\G$,
let $h=(h_1,\ldots, h_{k+1})$ be a triangular resolution
with
$$
	\ord_{z_j}h_j \le \mu,
	\quad
	1\le j\le k.
$$
Then 
$
	h_{k+1}\circ  \G
$
can be obtained
by holomorphic Kohn's procedures 
(P1) and (P2)
starting with the initial set consisting of components of 
$\psi$
and the ideal $I_k$,
where the number of procedures and the root order in (P2)
are effectively bounded 
by a function depending only on $(n,q,\mu)$.
\end{enumerate}
\et

The proof for each of the statements in (MP1), (MP2) and (MP3)
will be provided
respectively in \S\ref{jacs}
and Propositions~\ref{triangle} and \ref{claim}.
All three meta-procedures will be subsequently 
used one after another in \S\ref{6} to prove
the following explicit description of 
$q$-multipliers arising from our algorithm:

\bc\Label{k}
	For integers $n,q,\nu\ge 1$, initial system 
	$\psi_0=(\psi_{0,1}, \ldots, \psi_{0,n-q+1})$
	of $q$-multiplicity $ \le \nu$,
	and  $1\le k\le n-q+1$, there exist:
	\begin{enumerate}
	\item
		holomorphic coordinates $(z_1, \ldots, z_n)$ chosen among linear combinations of any given holomorphic coordinate system;
	\item
		systems 
		$\psi_k = (\psi_{k,k+1}, \ldots, \psi_{k,n-q+1})$
		chosen among generic linear combination of 
		$\psi_0$,
		and associated maps
		$$
			\G_k(z):= (z_1, \ldots, z_k, \psi_{k,k+1}(z), \ldots, \psi_{k,n-q+1}(z),z_{n-q+2},\ldots, z_n);
		$$
	\item
		systems of multipliers
		$f_k = (f_{k,1}, \ldots, f_{k,k})$
		obtained
		via effective meta-procedures applied
		to $(\psi_{k-1}, f_{k-1} )$ (where $f_0$ is empty);
	\item
		integer functions $\nu_{k,j}(n, q,\nu)>0$
		and
		decompositions
		of the form $f_{k,j} = Q_{k,j}\circ \G_k$, 
		$j=1, \ldots, k$,
		where 
		each $Q_{k,j}=Q_{k,j}(w_j, \ldots, w_n)$ is a holomorphic function 
		depending only on the last $n-j+1$ coordinates
		with $\ord_{w_j}Q_{k,j}\le \nu_{k,j}(n, q,\nu)$
		;
	\item 
		positive functions $\eps_{k,j}(n, q,\nu)>0$
		such that the order of subellipticity of each $f_{k,j}$ 
		 for $(0,q)$ forms
		is $\ge\eps_{k,j}(n,q, \nu)$.
	\end{enumerate}
\ec

Using Corollary~\ref{k}
for the largest $k$,
we prove
Theorem~\ref{p12}
by
 applying the 
 meta-procedure (MP3) from Theorem~\ref{meta}
one last time:

\bpf[Proof of Theorem~\ref{p12}]
Taking $k=n-q+1$ in Corollary~\ref{k},
we find 
$$\G_k(z)=(z_1,\ldots,z_k,z_{k+1},\ldots, z_n)=z,$$
whose Jacobian determinant $J=1$.
Hence, using $f_{k,j}=Q_{k,j}\circ\G_k$ provided by
Corollary~\ref{k}, we can 
apply (MP3) from Theorem~\ref{meta}
to
$$
\G=\G_k,
\quad
I_j = (f_{k,1},\ldots,f_{k,j})
\text{ for }
j\le k,
\quad
I_{k+1}=1 \subset I_k+ (J),
$$
and
$$
\quad
h_j=Q_j 
\text{ for } j\le k,
\quad
h_{k+1}=1,
$$
to conclude that 
$h_{k+1}\circ\G=1$ is a $q$-multiplier,
completing the proof.
\epf

\section{Preliminaries}
\subsection{Multiplicity and degree}
\Label{defs}
Denote by $\6O=\6O_{n,p}$
the ring of germs at a point $p$
of holomorphic functions in $\C^n$.
Since our considerations are 
for germs at a fixed point,
we shall assume $p=0$
unless specified otherwise.

Recall that an ideal
$\6I\subset\6O$
is of {\em finite type}
if $\dim\6O/\6I<\infty$,
or equivalently the (germ at $0$ of the) zero variety 
$\6V(\6I)$
is zero-dimensional at $0$.
In the latter case,
the classical {\em algebraic intersection multiplicity} of $\6I$
(see e.g.\ \cite[\S1.6, \S2.4]{Fu84})
is defined as
\beq
	\mult \6I := \dim \6O/\6I.
\eeq
Similarly, for a germ of holomorphic map
$\psi\colon(\C^n,0)\to (\C^n,0)$,
we have $\mult\psi := \mult (\psi)$,
where $(\psi)$ is the ideal generated
by the components of $\psi$,
and the quotient $\6O/(\psi)$ is the {\em local algebra of $\psi$}
(see e.g.\ \cite{AGV85}).
More generally (cf.~ \cite[\S2.4]{D93}), for every integer $0< q \le n$,
define the {\em $q$-multiplicity} by
\beq\Label{6I}
	\mult_q \6I := \min \dim  \6O/(\6I + (L_1, \ldots, L_{q-1})),
\eeq
where the minimum is taken over
sets of $q-1$  linear functions $L_j$ on $\C^n$.
The same minimum is achieved
when $L_j$ are germs
of holomorphic functions with linearly independent differentials,
as can be easily shown by a change of coordinates
linearlizing the functions.
In particular, the $q$-multiplicity of an ideal is a {\em biholomorphic invariant}.
In a similar vein, given a collection 
$\phi=(\phi_1, \ldots, \phi_{n-d})\in \6O_n$
of
$n-d$ function germs vanishing at $0$,
we write
\beq\Label{mult-d}
	\mult (\phi)
	=
	\mult (\phi_1, \ldots, \phi_{n-d}) 
	:= \min \dim  \6O/(\phi_1, \ldots, \phi_{n-d}, L_1, \ldots, L_d),
\eeq
where $L_j$ are as above.
That is, we will adopt the following convention:

\medskip

{\bf Convention.}
For every $1\le k\le n$ and a $k$-tuple 
of holomorphic function germs
$\phi_1,\ldots, \phi_k$,
their multiplicity 
$\mult(\phi_1,\ldots, \phi_k)$
is always assumed 
to be the $(n-k+1)$-multiplicity,
i.e.\ with $(n-k)$ generic linear functions
added to the ideal.
\medskip


Further recall that the {\em degree} 
$\deg(\psi)$
of a germ
(also called ``index" in \cite{AGV85}) 
of a finite holomorphic map
$\psi\colon(\C^n,0)\to (\C^n,0)$
is the minimum $m$ such that
$\psi$ restricts to a ramified $m$-sheeted covering
between neighborhoods of $0$ in $\C^n$.
%
Both integers are known to coincide (see e.g.\ \cite{ELT77, AGV85, D93}):

\begin{Thm}
[{\cite[\S4.3]{AGV85}}]
\Label{degree}	
Let $\psi\colon(\C^n,0)\to (\C^n,0)$ be germ of finite holomorphic map.
Then
$$
	\mult (\psi) = \deg \psi.
$$
\end{Thm}

\subsection{Basic properties of multiplicity}
The proofs of the following lemmas can be found in \cite{S10} and \cite{KZ20}.

\bl[Semicontinuity of multiplicity]\Label{semi}
Let $\psi_t \colon  (\C^m,0)\to (\C^k,0)$
be a continuous family of germs of holomorphic maps,
in the sense that all coefficients of the power series expansion of $\psi_t$
depend continuously on $t\in \R^m$.
Then $\mult(\psi_t)$ is upper semicontinuous in $t$.
\el

In the following we keep using the notation \eqref{mult-d}.

\bc\Label{ineq}
For every germs
$$
	(f,g)\colon (\C^{n+m}, 0) \to (\C^n, 0) \times (\C^m, 0),
$$
we have
$$
	\mult f \le \mult (f,g).
$$
\ec



\bd[D'Angelo, \cite{D82}]
Let $S\subset \mathcal O_{n,0}$ be a subset of germs of holomorphic functions.
\ben
\item 
The D'Angelo $1$-type of $S$ is
$$
	\Delta^1 (S):=\sup_\g \inf_{f\in S} \frac{\ord f\circ\g }{\ord \g},
$$
where $\ord$ denotes the vanishing order, and
the supremum is taken over all nonzero germs of holomorphic maps
$\g\colon(\C,0)\to(\C^n,0)$.
\item
the D'Angelo $q$-type of $S$ for $q\ge1$ is
$$
	\Delta^q (S):=\inf _L ~\Delta^1(S\cup L),
$$
where the infimum is taken over all sets $L$ of $(q-1)$ complex linear functions.
\een
\ed

Let $\Omega$ be a domain defined locally by
\beq\Label{om1}
	\Re(z_{n+1})+\sum_{j=1}^N |F_j(z_1,\ldots, z_n)|^2<0,
\eeq
where $F_1,\ldots,F_N$ are holomorphic functions
in a neighborhood of $0$.
By $q$-{\emph type} $\Delta^q(\Omega)$ of \eqref{om1} at $0$ 
we mean twice the $q$-type of $F_1,\ldots, F_N$.
Let $p$ be the smallest integer such that
$$
	|z|^p\leq A\sum _j |F_j|+|L(z)|
$$
holds for some linear map $L:\mathbb C^n\to \mathbb C^{q-1}$. By following the argument in (I.2) of \cite{S10}, we can show that 
$$
	2p=\Delta^q(\Omega)\leq \mult_q \mathcal I(F_1,\ldots,F_N)=:s.
$$
Furthermore, for the smallest integer $r$ satisfying
$$
	\mathfrak m^r\subset \mathcal I(F_1,\ldots, F_N,  L_1,\ldots, L_{q-1})
$$
for some linear functions $L_1,\ldots,L_{q-1}$, where $\mathfrak m$ is the maximal ideal in $\mathcal O_{n,0}$,
we obtain
$$
	r\leq s\leq (n+r-1)!/n!( r-1)!.
$$

An important ingredient is the following 
consequence from Siu's lemma
on selection of linear combinations
of holomorphic functions
 for effective multiplicity
 \cite[(III.3)]{S10}
 combined with effective comparison
 of the invariants of holomorphic map germs
 \cite[(I.3-4)]{S10}
 (see also \cite[\S2.2]{D93}):
\bl[$q$-type version of Siu's lemma on effective mixed multiplicity]\Label{siu}
Let $0\le j \le n-q$ and
 $f_1,\ldots, f_j, F_1, \ldots, F_N$
be holomorphic function germs in $\6O_{n,0}$ 
such that
$$
	\mu:= \mult(f_1, \ldots, f_j) <\infty,
	\quad
	\nu:=\mult_q(F_1,\ldots, F_N) <\infty.
$$
Then
$$
	\mult (f_1, \ldots, f_j, G_{j+1}, \ldots, G_{\ell}) 
	\le
	\mu  \nu^{\ell-j},\quad \ell\leq n-q+1
$$
holds for generic linear combinations
$G_{j+1}, \ldots, G_{\ell}$ of $F_k$'s.
\el


We shall also need the following 
lemma proved in \cite[\S3]{KZ20}:
\begin{Lem}[Effective Nullstellensatz, \cite{KZ20}]\Label{radical}
Let $\phi_1,\ldots,\phi_k,f\in \mathcal{O}_{n,0}$
satisfy
\begin{equation*}
	\mu := \mult(\phi_1,\ldots,\phi_k)<\infty,
		\quad
	f\in \sqrt{(\phi_1,\ldots,\phi_k)}.
\end{equation*}
Then
$$
	f^{n\mu }\in (\phi_1,\ldots,\phi_k).
$$
\end{Lem}

\section{Multiplicity estimates for Jacobian determinants}\Label{jacs}


The meta-procedure (MP1) in Theorem~\ref{meta} is the special case of the following proposition where one can put
$d=n-q-k$ and 
$\psi=(\psi_{k+1},\ldots,\psi_{n-q+1}, z_{n-q+2},\ldots,z_n)$ assuming that
$$
	\mult_q(f_1,\ldots,f_{k},\psi_{k+1},\ldots,\psi_{n-q+1})=\mult(f_1,\ldots,f_{k},\psi_{k+1},\ldots,\psi_{n-q+1},z_{n-q+2},\ldots,z_n)<\infty
$$
after a suitable linear coordinate change of $(z_1,\ldots,z_n)$. The proof of the proposition is given in \cite{KZ20}. In what follows,
we use the convention that
$\mult (f)=1$ if $f$ has $0$ components.

\bp[\cite{KZ20}]
\Label{multiplicity-J}
Let
$(f,\psi)\colon (\C^n,0)\to (\C^{n-d}\times\C^d,0)$,
$1\le d\le n$,
be a holomorphic map germ satisfying
$$
	\nu:=\mult (f) <\infty,
	\quad 	
	\mu:= \mult (f, \psi)  <\infty.
$$
Then after a linear change of 
coordinates 
$(z_1,\ldots,z_{n})$
and another linear coordinate change in $\C^d$,
the partial Jacobian determinant
\beq\Label{J-part}
	J:=\frac{
		\partial(\psi_1,\ldots,\psi_{d})
	}{
		\partial(z_1,\ldots,z_{d})
	}
\eeq
satisfies
\beq\Label{fjpsi}
    \mult(f, J)\le d \nu\mu,
    \quad
    \mult(f, J, \psi_2,\ldots,\psi_d)\le d \nu\mu^{d},
\eeq
where $\psi_j$ is the $j$-th component of $\psi$ in the new coordinates.
\ep

\section{Existence of effective triangular resolutions}

The following is a more precise version
of the meta-procedure (MP2) in Theorem~\ref{meta}.
We shall denote by $\ord_{w_j} h$
the vanishing order at $0$ of $h(0,\ldots,0,w_j,0,\ldots,0)$
(where all variables are zero except $w_j$).

\bp\Label{triangle}
Let $1\le k\le n-q+1$, 
 $f_1, \ldots, f_k, \phi_1,\ldots, \phi_{n-q+1} \in \6O_{n,0}$, 
satisfy
\beq\nonumber
	\mu_j:=
	\mult_q (f_1, \ldots, f_j, \phi_{j+1}, \ldots, \phi_{n-q+1}) <\infty,
	\quad
	1\le j\le k.
\eeq
Let $\mathcal I$ be the filtration of ideals 
$$I_j:=(f_1,\ldots,f_j),~1\leq j\leq k.$$
Then
there exist a germ of a holomorphic map 
$$
\Gamma_\phi=(\phi_1,\ldots,\phi_{n-q+1},L_{n-q+2},\ldots,L_n)
$$ 
and a triangular resolution $h=(h_1,\ldots, h_k)$ of $(\Gamma_\phi, \mathcal I)$
such that 
\beq\Label{h-mult}
	\ord_{w_j} h_j \le 
	n\cdot \mu_j \cdot \mult(f_1,\ldots, f_j) 
	,
	\quad
	1\le j\le k.
\eeq
Furthermore, each $h_j(w_j, \ldots, w_n)$ can be chosen
as Weierstrass polynomial in $w_j$.
\ep

\bpf
Since 
$$
	\mult(\psi_1,\ldots,\psi_k)=\mult(\psi_1,\ldots,\psi_k,L_{k+1},\ldots,L_n),
$$
for generic choice of $n-k$ linear functions $L_j$, 
we can choose a set of linear functions $L_{n-q+2},\ldots,L_n$ such that
\beq\Label{muj}
	\mult(f_1,\ldots,f_j,\phi_{j+1},\ldots,\phi_{n-q+1},L_{n-q+2},\ldots,L_n)=\mu_j,\quad \text{ for all } j.
\eeq 
Let 
$$
\Gamma_\phi:=(\phi_1,\ldots, \phi_{n-q+1}, L_{n-q+2},\ldots, L_n).
$$ 
Consider 
the coordinate projections
$$
    \pi_j(w_1,\ldots, w_n)
    = 
    (w_{j}, \ldots, w_n)
    \in \C^{n-j+1}
    ,
    \quad
    1\le j \le k,
$$
and let
$$
	W_j := \6V \left(f_1,\ldots,f_j\right),
	\quad
	\widetilde W_j
	:=
	(\pi_{j} \circ \G_\phi)(W_j)
	\subset\mathbb{C}^{n-j+1},\quad 1\leq j\leq k,
$$
where $\6V$ is the zero variety.
Then $W_j$ is of codimension $\ge k$ in $\cn$.
In fact, counting preimages and using \eqref{muj}, we 
conclude that
$\widetilde W_{j}\subset \mathbb{C}^{n-j+1}$ is a proper subvariety of codimension $1$
and
$$
	\pi_{j+1} |_{\widetilde W_j} \colon \widetilde W_j \to \mathbb{C}^{n-j}
$$
is a finite 
holomorphic 
map germ
of degree
$\le \mu_j$.
Then there exist  Weierstrass polynomials 
$Q_j(w_j, \ldots, w_n)$, $j=1,\ldots,k$, 
satisfying
$$
	Q_j=
	w_{j}^{\nu_j}+\sum_{\ell<\nu_j}b_{j,\ell}(w_{j+1},\ldots,w_{n})w_{j}^\ell,
	\quad
	Q_j|_{\2W_j} =0,
	\quad
	\nu_j
	=
	\ord_{w_j} Q_j
	\le \mu_j.
%
$$
Furthermore, Lemma~\ref{radical} implies
$$
	h_j\circ\G_\phi
	\in(f_1,\ldots,f_j),
	\quad
	h_j:= Q_j^{\l_j},
$$
for suitable $\l_j\in \N_{\ge1}$ satisfying
$$
	\l_j \le n\cdot \mult(f_1,\ldots, f_j).
$$
Then $(h_1,\ldots,h_k)$ is a triangular resolution
satisfying \eqref{h-mult}
 as desired.
\end{proof}

\section{Effective Kohn's procedures for triangular resolutions}

The following is a more precise version
of the meta-procedure (MP3) in Theorem~\ref{meta}:

\bp\Label{claim}
Let $1\le k \le n-q$ and $(Q_1\circ \Gamma,\ldots,Q_{k+1}\circ\Gamma)$ 
be a triangular resolution 
of $(\G,\6I)$,
where 
$\G\colon (\cn,0)\to(\C^n,0)$
is a holomorphic map germ of the form 
\beq
\Label{GQ}
\Gamma=(\psi_1,\ldots,\psi_{n-q+1},z_{n-q+2},\ldots,z_n),
\eeq
and $\6I$
a filtration of ideals $I_1\subset \ldots \subset I_{k+1}\subset\6O_{n,0}$.
Assume
\beq\Label{mus}
	\mu_j = \ord_{w_j} Q_j <\infty,
	\quad
	1\le j\le k,
\eeq
and
\beq\Label{det-inc}
	I_{k+1}\subset  I_k+(J),
\eeq
where $J$ is the Jacobian determinant of $\G$.

Then 
$Q_{k+1}\circ \G$ can be obtained
by applying holomorphic Kohn's procedures 
(P1) and (P2)
to $(\Gamma, I_k)$
and
each procedure (P1) and (P2) is applied $\mu_1\cdots \mu_k$ number of times with the root order in (P2) being $\le k+1$.
In particular, if $I_k$ consists of multipliers of order $\ge \eps$,
then $Q_{k+1}\circ \G$ is a multiplier of order $\geq (2k+2)^{-\mu_1\cdots\mu_k} \eps$.
\ep

\begin{proof}
Since $\mu_j<\infty$ for $j\le k$, 
multiplying by invertible holomorphic functions,
we may assume that 
$$
	Q_j=
	w_{j}^{\mu_j}
	+\sum_{\ell<\mu_j}b_{j,\ell}(w_{j+1},\ldots,w_{n})w_{j}^\ell,
	\quad
	1\le j\le k,
$$
are Weierstrass polynomials satisfying
$$
	f_j:= Q_j\circ \G\in I_j.
$$
In addition, \eqref{det-inc} implies
\beq\Label{f-inc}
	f_{k+1} 
	:=
	Q_{k+1}\circ\G \in I_k + (J).
\eeq

For simplicity of notation,
we use the remaining indices to
denote the coordinate functions
in \eqref{GQ}, i.e.\
$$
\Gamma=(\psi_1,\ldots,\psi_{n-q+1},z_{n-q+2},\ldots,z_n)
=
(\psi_1,\ldots,\psi_{n-q+1},\psi_{n-q+2},\ldots,\psi_n)
.
$$
Since
$$
	(f_1, \ldots, f_k,  \psi_{k+1},\ldots,\psi_{n} ) 
	= \Phi \circ \G,
$$
where
\beq\Label{P}
	\Phi(w) 
		:= 
	(Q_1(w), \ldots, Q_k(w), w_{k+1}, \ldots, w_n),
\eeq
we obtain the Jacobian determinants
$$
	\frac{\partial( f_1,\ldots, f_k,
	 \psi_{k+1}, \ldots,  \psi_{n-q+1})}
	{\d (z_1,\ldots, z_{n-q+1})}
	=
	\frac{\d(f_1,\ldots, f_k, \psi_{k+1},\ldots,\psi_{n})}
	{\d (z_1,\ldots, z_{n})}
	=: J_{(1,\ldots,1)}.
$$
For $L=(\ell_1,\ldots,\ell_k)\in \mathbb N^k$,
 define $A_L\in \6O_{n,0}$ by
$$
	A_L(w):=\partial_{w_1}^{\ell_1}Q_1(w) \cdots\partial_{w_k}^{\ell_k}Q_k(w).
$$
%
%
Then the Jacobian factors as 
$$
	J_{(1,\ldots,1)}
	=
	(A_{(1,\ldots,1)}\circ \G) J
$$
and hence by \eqref{f-inc},
$$
	 (A_{(1,\ldots,1)} \circ \Gamma)f_{k+1} 
	 \in I_k + (J_{(1,\ldots,1)})
$$
is obtained by applying Kohn's procedure (P1)
to $I_k$.

Now for $B=(B_1,\ldots,B_k)$, where 
$$
	B_j:=(A_{L_j} \circ \Gamma)f_{k+1} 
	= 
	(A_{L_j} Q_{k+1}) \circ \Gamma
$$
or
$$
	B_j:=Q_j \circ \Gamma =f_j
$$
with each $B_j$ in both cases
are obtained by applying Kohn's procedures,
we obtain
$$
	\frac
	{ \partial (B_1,\ldots, B_k, \psi_{k+1},\ldots,  \psi_{n-q+1}) }
	{\d(z_1,\ldots, z_{n-q+1})}
	=
	\frac{\d(B_1,\ldots, B_k, \psi_{k+1},\ldots,\psi_{n})}
	{\d(z_1,\ldots, z_{n})}
	=: J_L.
$$
In view of our assumption that each $Q_j$, $j\le k$, 
is a Weierstrass polynomial in $w_j$ of degree $\mu_j$,
each top derivative $\d_{w_j}^{\mu_j}Q_j$ is constant
and hence $B_j$ only depends on $(w_j, \ldots, w_n)$.
Then using factorization of the Jacobian determinant
 and the triangular property of $B_j$'s,  
we obtain
$$
	J_L=
	c\left(
		\big(
		(\partial_{w_1}^{\ell_1}Q_1)^{m_1}
            \cdots(\partial_{w_k}^{\ell_k}Q_k)^{m_k}
		Q_{k+1}^{m_{k+1}}
		\big)
		  \circ
        \Gamma
        	\right)
	J
$$
for some constant $c\ne0$ and integers $m_j$,~$j=1,\ldots, k+1$
and hence by \eqref{f-inc},
$$
	 \left(
	 	(A_L Q_{k+1}) \circ \Gamma
	\right)^{m_{k+1}+1}
	 \in I_k + (J_L).
$$
Then  $(A_L Q_{k+1})\circ\G$ is obtained by the Kohn's procedure (P2)
with root order $\le m_{k+1}+1$,
and by using the lexicographic order for $L=(\ell_1,\ldots,\ell_k)$ as in the proof \cite{KZ20}, we can complete the proof.
\end{proof}

%
%
%
%

\section{Proof of Corollary \ref{k}}
\Label{6}

We will use the induction on $k$.
For the case $k=1$, take
 $\psi_0=(\psi_{0,1},\ldots,\psi_{0,n-q+1})$ to be a (generic)
linear combinations
of $F_j$'s
such that $\mult_q(\psi_0)$ is effectively bounded 
and assume that
$$	
	\mult_q(\psi_0)=\mult (\psi_0, z_{n-q+2},\ldots,z_n)
$$
after a linear coordinate change of $\mathbb C^n$. Such $\psi_0$ exists by Lemma \ref{siu}.

Now suppose that the statement of the corollary holds for $k-1$. 
Applying Lemma~\ref{siu} and (MP1),
we obtain (generic)
linear combinations
$\psi_k=(\psi_{k,k+1},\ldots,\psi_{k,n-q+1})$
 of $\psi_{k-1, j}$'s
such that $
	\mult_q(z_1,\ldots, z_k,\psi_k)$
and $\mult_q(f_{k-1},J, \psi_k)$
are effectively bounded, where 
$$
	J:= 
		\frac{
			\partial(\psi_{k-1,k},\ldots,\psi_{k-1,n-q+1})
		}{
			\partial (z_{k},\ldots,z_{n-q+1})
		} .
$$

Next apply (MP2)
for the map germ
$$
	\G(z):=(z_1,\ldots, z_{k-1},\psi_{k-1}(z),z_{n-q+2},\ldots,z_n)
$$ 
and the filtration of ideals 
$$
	I_j:=(f_{k-1,1},\ldots, f_{k-1,j}), 
		\; 
	1\le j\le k-1,
$$
to obtain a triangular resolution $(h_1,\ldots, h_{k-1})$ 
such that 
$\ord_{w_{j}} h_{j}$ is effectively bounded.

Finally, apply (MP3) for 
$$
	(\phi,\psi)=(z_1,\ldots,z_{k-1}, \psi_{k-1})
$$ 
and a filtration $\mathcal I $ of ideals 
$$
	\tilde I_j=(h_1\circ \Gamma,\ldots, h_j\circ\Gamma),\quad
	j=1,\ldots,k-1
$$
and 
$$
	\tilde I_k=\tilde I_{k-1}+(\tilde J)\subset \tilde I_{k-1}+(J),
$$
%
%
where
$$
	\tilde J=
	\frac
	{\partial( (h_1\circ\Gamma),\ldots,(h_{k-1}\circ\Gamma),
	\psi_{k-1,k},\ldots,\psi_{k-1,n-q+1})}
	{\d(z_1,\ldots,z_{n-q+1})}
.
$$
Then we obtain a new set of 
multipliers
$$
	f_k=(f_{k,1},\ldots,f_{k,k})
$$
given by the triangular resolution $(h_1,\ldots,h_k)$ of $(\Gamma, \mathcal I)$
together with a set of premultipliers
$$
	\psi_k=(\psi_{k,k+1},\ldots,\psi_{k,n-q+1})
$$
that satisfy the condition of the corollary,
completing the proof.

\end{document}